\documentclass{amsart}
\usepackage{graphicx}
\usepackage{amsfonts}
\usepackage{amsmath}
\usepackage{amssymb}
\newtheorem{thm}{Theorem}[section]
\newtheorem{lem}[thm]{Lemma}
\newtheorem{cor}[thm]{Corollary}

\newtheorem{obs}[thm]{Observation}
\theoremstyle{definition}
\newtheorem{defn}[thm]{Definition}

\theoremstyle{remark}
\newtheorem{remark}[thm]{Remark}
\numberwithin{equation}{section}
\DeclareMathOperator{\dd}{\Lambda}
\DeclareMathOperator{\Q}{Qsym}
\DeclareMathOperator{\BQ}{BQsym}

\DeclareMathOperator{\Des}{Des}

\DeclareMathOperator{\Pe}{Pk}

\DeclareMathOperator{\lPe}{Pk^{(\ell)}}

\title[Enriched $P$-partitions and peak algebras]{Enriched $P$-partitions and peak algebras\\(extended abstract)}

\author[T. K. Petersen]{T. Kyle Petersen}
\address{Department of Mathematics, Brandeis University, Waltham, MA, USA, 02454}
\email{tkpeters@brandeis.edu}
\urladdr{http://people.brandeis.edu/\~{}tkpeters}


\begin{document}
\begin{abstract}
We generalize Stembridge's enriched $P$-partitions and use this theory to outline the structure of peak algebras for the symmetric group and the hyperoctahedral group. Whereas Stembridge's enriched $P$-partitions are related to quasisymmetric functions (the coalgebra dual to Solomon's type A descent algebra), our generalized enriched $P$-partitions are related to type B quasisymmetric functions (the coalgebra dual to Solomon's type B descent algebra). Using these functions, we explore three different peak algebras: the ``interior" and ``left" peak algebras of type A, and a new type B peak algebra. Our results specialize to results for commutative peak algebras as well.
\end{abstract}

\maketitle

\section{Introduction}

Much attention has been given to the so-called \emph{descent algebras}; see \cite{BergeronBergeron, BergeronBergeron2, Bergeron, Cellini, Cellini2, Cellini3, Fulman, GarsiaReutenauer, Loday, Mahajan, Petersen, Solomon}. Here we add a chapter to the story of the more recently introduced \emph{peak algebras}. Our approach expands on the one taken in \cite{Chow, Gessel} and \cite{Petersen}, where descents were studied using Richard Stanley's $P$-partitions, or modified versions thereof. This paper is a condensed version of \cite{Petersen2}, which contains several results not mentioned here, as well as any omitted proofs.

Generically, a \emph{peak} of a permutation $\pi \in \mathfrak{S}_n$ is a position $i$ such that $\pi(i-1) < \pi(i) > \pi(i+1)$. The only difference between the various types of peak sets we will study is the values of $i$ that we allow. The \emph{interior peak set} and the \emph{left peak set} are, respectively:
\begin{align*}
\Pe(\pi) & := \{ i \in [2,n-1] \,|\, \pi(i-1) < \pi(i) > \pi(i+1) \} \\
\lPe(\pi) & := \{ i \in [1,n-1] \,|\, \pi(i-1) < \pi(i) > \pi(i+1) \},
\end{align*}
where we take $\pi(0) = 0$. For example, the permutation $\pi = (2,1,4,3,5)$ has $\Pe(\pi) = \{3\}$, $\lPe(\pi) = \{1,3\}$. We will also study the peak set of signed permutations $\pi \in \mathfrak{B}_n$, defined by \[ \Pe_B(\pi) := \{ i \in [0,n-1] \,|\, \pi(i-1) < \pi(i) > \pi(i+1) \},\] where $\pi(0) = 0$ and we say there is a peak in position $0$ if $\pi(1) < 0$. For example, if $\pi = (-2,3,4,-5,1)$, then $\Pe_B(\pi) = \{ 0, 3 \}$. One can study the suitably defined \emph{right} and \emph{exterior} peaks as well, but the algebraic implications are more limited. See sections \ref{sec:spec} and \ref{sec:neg}, and also \cite{Petersen2}.

The study of algebras related to peaks began with John Stembridge's paper \cite{Stembridge} on \emph{enriched} $P$-partitions, followed by others, including \cite{AguiarBergeronNyman, AguiarNymanOrellana, BergeronHivertThibon, BergeronHohlweg, BergeronMykytiukSottileWilligenburg, BilleraHsiaoWilligenburg, KrobThibon, Schocker}. While \cite{Stembridge} explores ``the algebra of peaks" related to quasisymmetric functions, it does not use enriched $P$-partitions for the study of subalgebras of the group algebra $\mathbb{Z}[\mathfrak{S}_n]$ as we will here, and the only notion of peak that it uses is that of an interior peak. Kathryn Nyman \cite{Nyman} built on \cite{Stembridge} to show that there is a subalgebra of the group algebra of the symmetric group, akin to Solomon's descent algebra \cite{Solomon}, formed by the linear span of \[ v_I := \sum_{\substack{\pi \in \mathfrak{S}_n \\ \Pe(\pi) = I }} \pi,\] which we call the \emph{interior peak algebra}, denoted $\mathfrak{P}_n$. Later, without the use of enriched $P$-partitions, Marcelo Aguiar, Nantel Bergeron, and Nyman \cite{AguiarBergeronNyman} showed that left peaks also give a subalgebra in this sense. We will denote the linear span of sums of permutations with the same set of left peaks by $\mathfrak{P}^{(\ell)}_n$. In \cite{AguiarBergeronNyman}, the authors also examined commutative subalgebras of the peak algebras---the ``Eulerian" peak algebras formed by sums of permutations with the same number of peaks. One goal of this work is to derive some of the results of \cite{AguiarBergeronNyman} as a natural application of enriched $P$-partitions. In doing so, we are led to the type B enriched $P$-partitions and to the type B peak algebra, $\mathfrak{P}_{B,n}$.

The link between peak algebras and enriched $P$-partitions is through quasisymmetric generating functions. Let $\Q := \bigoplus_{n\geq 0} \Q_n$ denote the space of quasisymmetric functions, where $\Q_n$ denotes the quasisymmetric functions homogeneous of degree $n$. Ira Gessel \cite{Gessel} showed how generating functions for ordinary $P$-partitions give a natural basis for $\Q$, and moreover, he defined a coproduct on $\Q_n$ that makes it the coalgebra dual to Solomon's descent algebra for the Coxeter group of type $A_{n-1}$. Stembridge \cite{Stembridge} defined generating functions for enriched $P$-partitions that form a subring of the ring of quasisymmetric functions, called the \emph{peak functions}. Let $\mathbf{\Pi} := \bigoplus_{n\geq 0} \mathbf{\Pi}_n$ denote the space of peak functions, with $\mathbf{\Pi}_n$ the $n$-th graded component. We will use an approach similar to Gessel's to give a coproduct on $\mathbf{\Pi}_n$ that makes it dual to Nyman's interior peak algebra.

Just as Stembridge's enriched $P$-partitions connect with quasisymmetric functions (the coalgebra dual to Solomon's type A descent algebra), the new types of enriched $P$-partitions we present here connect to the type B quasisymmetric functions, $\BQ:= \bigoplus_{n \geq 0} \BQ_n$ (the coalgebra dual to Solomon's type B descent algebra), as defined by Chak-On Chow \cite{Chow} using type B $P$-partitions. We will define the \emph{type B peak functions} $\mathbf{\Pi}_{B} := \bigoplus_{n \geq 0} \mathbf{\Pi}_{B,n}$ and the \emph{left peak functions} $\mathbf{\Pi}^{(\ell)} := \bigoplus_{n \geq 0} \mathbf{\Pi}^{(\ell)}_n$, and give a natural coproduct that makes $\mathbf{\Pi}^{(\ell)}_n$ dual to $\mathfrak{P}^{(\ell)}_n$ and $\mathbf{\Pi}_{B,n}$ dual to $\mathfrak{P}_{B,n}$.

\begin{remark}
It is known that the quasisymmetric functions form a Hopf algebra, and Stembridge's peak functions are a Hopf subalgebra \cite{BergeronMykytiukSottileWilligenburg}. A natural question is whether the type B quasisymmetric functions form a Hopf algebra, and they do. As of this writing, it is known that the left peak functions do \emph{not} form a Hopf subalgebra, but an unresolved question is whether type B peak functions form a Hopf subalgebra. This topic is part of ongoing work.
\end{remark}

\section{Enriched $P$-partitions}\label{sec:epp}

\begin{figure} [h]
\centering
\includegraphics{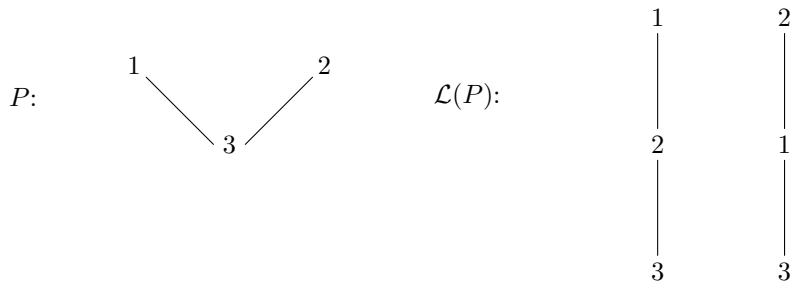}
\caption{Linear extensions of a poset $P$.\label{fig:hasse}}
\end{figure}

The ``$P$" in $P$-partition stands for a partially ordered set, or poset. For our purposes, we assume that all posets $P$, with partial order $<_{P}$, are finite. And unless otherwise noted, if $|P| = n$, then the elements of $P$ are labeled distinctly with the numbers $1,2,\ldots,n$. We will sometimes describe a poset by its Hasse diagram, as in Figure \ref{fig:hasse}. We can think of any permutation $\pi \in \mathfrak{S}_{n}$ as a poset with the total order $\pi(s) <_{\pi} \pi(s+1)$.

For a poset $P$ with $n$ elements, let $\mathcal{L}(P)$ denote its Jordan-H\"{o}lder set: the set of all
permutations of $[n]$ which extend $P$ to a total order. This set is also called the set
of ``linear extensions" of $P$. For example let $P$ be the poset
defined by $1 >_{P} 3 <_{P} 2$. In linearizing $P$ we form a
total order by retaining all the relations of $P$ but introducing
new relations so that any element is comparable to any other. In
this case, 1 and 2 are not comparable, so we have exactly two ways
of linearizing $P$: $3 < 2 < 1$ or $3 < 1 < 2$. These correspond
to the permutations $(3,2,1)$ and $(3,1,2)$. Let us make the
following observation.
\begin{obs}\label{ob1} A permutation $\pi$ is in $\mathcal{L}(P)$ if and only if $i<_{P}
j$ implies $\pi^{-1}(i) < \pi^{-1}(j)$.
\end{obs}
In other words, if $i$ is ``below" $j$ in the Hasse diagram of the poset $P$, it must be below $j$ in any linear extension of the poset.

We now introduce the basic theory of enriched $P$-partitions, building on Stembridge's work \cite{Stembridge}. To begin, Stembridge defines $\mathbb{P}'$ to be the set of nonzero integers with the following total order: \[ -1 < 1 < -2 < 2 < -3 < 3 < \cdots\] We will have use for this set, but we view it as a subset of a similar set. Define $\mathbb{P}^{(\ell)}$ to be the integers with the following total order: \[0 < -1 < 1 < -2 < 2 < -3 < 3 < \cdots\] Then $\mathbb{P}'$ is simply the set of all $i \in \mathbb{P}^{(\ell)}$, $i > 0$. In general, for any countable totally ordered set $S = \{s_1, s_2, \ldots\}$ we define $S^{(\ell)}$ to be the set \[\{s_0, -s_1, s_1, -s_2, s_2, \ldots\},\] with total order \[s_0 < -s_1 < s_1 < -s_2 < s_2 < \cdots \] (so we can think of $S^{(\ell)}$ as two interwoven copies of $S$ along with a zero element) and define $S'$ to be the set $\{ s \in S^{(\ell)} \,|\, s > s_0 \}$. For any $s_i \in \{s_0\} \cup S$, we say $s_i \geq 0$, or $s_i$ is \emph{nonnegative}. On the other hand, if $i\neq 0$ we say $-s_i < 0$ and $-s_i$ is \emph{negative}. The absolute value removes any minus signs: $|\pm s| = s$ for any $s \in \{s_0\} \cup S$.

For $i$ and $j$ in $S^{(\ell)}$, we write $s\leq^{+} t$ to mean either $s < t$ in $S^{(\ell)}$, or $s = t \geq 0$. Similarly we define $s \leq^{-} t$ to mean either $s < t$ in $S^{(\ell)}$, or $s = t < 0$. For example, on $\mathbb{P}^{(\ell)}$, we have $\{ s \, | \, s \leq^{+} 3 \} = \{\,0, \pm 1, \pm 2, \pm 3 \,\}$, $\{\, s \, | \, s \leq^{-} 3 \,\} = \{\, 0, \pm 1, \pm 2, -3\,\} = \{\, s \,|\, s \leq^{-} -3 \,\}$, $\{ \,s \,| \,0 \leq^{+} s \leq^{+} 2 \,\} = \{\, 0, \pm 1, \pm 2 \,\}$ and $\{\, s \,| \,0 \leq^{-} s \leq^{+} 2 \,\} = \{ \,\pm 1, \pm 2\, \}$.

\begin{defn}[Enriched $P$-partition]\label{def:epp}
An \emph{enriched $P$-partition} (resp. \emph{left enriched $P$-partition}) is an order-preserving map $f: P \to S'$ (resp. $S^{(\ell)}$) such that for all $i <_{P} j$ in $P$,
\begin{enumerate}
\item $f(i) \leq^{+} f(j)$ only if $ i < j$ in $\mathbb{Z}$,
\item $f(i) \leq^{-} f(j)$ only if $ i > j$ in $\mathbb{Z}$.
\end{enumerate}
\end{defn}

It is helpful to remember that Stembridge's enriched $P$-partitions are the nonzero left enriched $P$-partitions. We let $\mathcal{E}(P;S)$ denote the set of all enriched $P$-partitions $f: P \to S'$; $\mathcal{E}^{(\ell)}(P;S)$ denotes the set of left enriched $P$-partitions $f: P \to S^{(\ell)}$. If $S$ is irrelevant or understood, we simply write $\mathcal{E}(P)$ or $\mathcal{E}^{(\ell)}(P)$. For example, if our poset is $1 >_{P} 3 <_{P} 2$, then \[ \mathcal{E}^{(\ell)}(P) = \{ f: P \to S^{(\ell)} \,|\, f(1) \geq^- f(3) \leq^- 2 \}, \] which we can see actually splits into the two following disjoint subsets: \[ \{ f(3) \leq^- f(1) \leq^+ f(2) \} \sqcup \{ f(3) \leq^- f(2) \leq^- f(1) \} = \mathcal{E}^{(\ell)}( 312 ) \sqcup \mathcal{E}^{(\ell)}( 321).\]

This example leads us to the following, which, by analogy with a similar result for ordinary $P$-partitions, is referred to as the fundamental lemma of enriched $P$-partitions. It follows by induction on the number of incomparable pairs of elements in the poset.

\begin{lem}
\label{lem:FLEPP}
For any poset $P$, the set of all (left) enriched $P$-partitions is the
disjoint union of all (left) enriched $\pi$-partitions for linear
extensions $\pi$ of $P$. Equivalently,
\begin{align*}
\mathcal{E}(P) & = \coprod_{\pi \in \mathcal{L}(P)} \mathcal{E}(\pi),\\
\mathcal{E}^{(\ell)}(P) & = \coprod_{\pi \in \mathcal{L}(P)} \mathcal{E}^{(\ell)}(\pi).
\end{align*}
\end{lem}

Therefore when studying enriched $P$-partitions it is enough to consider the case where $P$ is a totally ordered chain, i.e., a permutation $\pi$. It is easy to describe the set of all enriched $\pi$-partitions in terms of descent sets. For any $\pi \in \mathfrak{S}_n$ we have
\begin{equation}
\label{eq:epp}
\begin{aligned}
\mathcal{E}(\pi) = \{\, f: [n] \to S' & \mid f(\pi(1)) \leq f(\pi(2)) \leq \cdots \leq f(\pi(n)), \\
&  i \notin \Des(\pi) \Rightarrow f(\pi(i)) \leq^{+} f(\pi(i+1))  \\
&  i \in \Des(\pi) \Rightarrow f(\pi(i)) \leq^{-} f(\pi(i+1))  \, \},
\end{aligned}
\end{equation}
and the analogous description for $\mathcal{E}^{(\ell)}(\pi)$ where we replace $S'$ with $S^{(\ell)}$.

From \eqref{eq:epp} it is clear that enriched $\pi$-partitions depend on the descent set of $\pi$. The connection to peaks is less obvious. In section \ref{sec:qsym} we will establish this link, and also show how left enriched $\pi$-partitions are related to left peaks. First, we present our main theorem.

\begin{figure} [h]
\centering
\includegraphics[scale = .7]{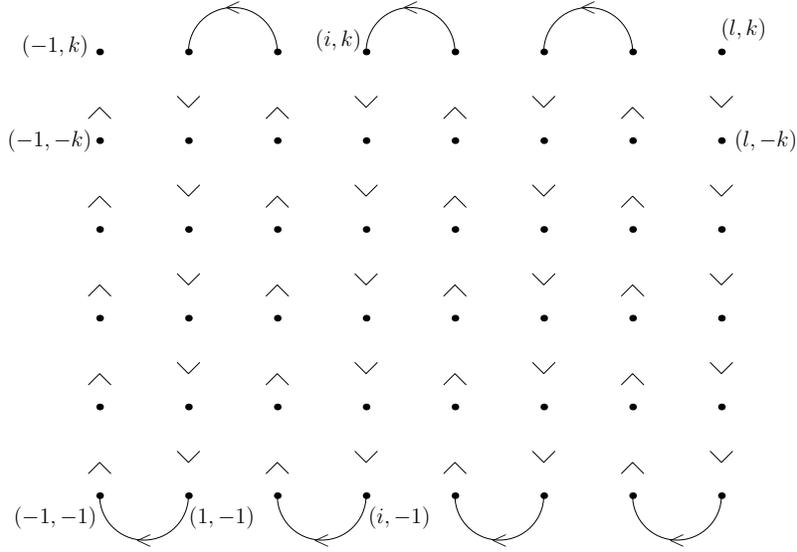}
\caption{The up-down order for $\mathbb{P}^{(\ell)} \times \mathbb{P}^{(\ell)}$.\label{fig:updown}}
\end{figure}

Let $S$ and $T$ be any two countable totally ordered sets, and let $S' \times T' = \{ (s,t) \,|\, s \in S', t \in T'\}$ be the cartesian product of $S'$ and $T'$ with the \emph{up-down order} defined as follows: $(s,t) < (u,v)$ if and only if
\begin{enumerate}
\item $s < u$, or

\item $s = u > 0$ and $t < v$, or

\item $s = u < 0$ and $t > v$.
\end{enumerate}
In other words, we read up the nonnegative columns, down the negative ones. Here we write $(s,t) \leq^{+} (u,v)$ in one of three cases: if $s<u$, or if $s = u > 0$ and $t \leq^{+} v$, or if $s = u < 0$ and $t \geq^{-} v$. Similarly, $(s,t) \leq^{-} (u,v)$ if $s < u$, or if $s=u > 0$ and $t \leq^{-} v$, or if $s = u < 0$ and $t \geq^{+} v$. We define $S^{(\ell)}\times T^{(\ell)}$ in the same way. See Figure \ref{fig:updown}.

\begin{thm}\label{thm:bijA}
We have the following bijections:
\begin{align}
\mathcal{E}(\pi;S\times T) & \longleftrightarrow \coprod_{\sigma\tau = \pi} \mathcal{E}(\tau; S) \times \mathcal{E}(\sigma; T) \label{eq:bijA1}\\
\mathcal{E}^{(\ell)}(\pi;S\times T) & \longleftrightarrow \coprod_{\sigma\tau = \pi} \mathcal{E}^{(\ell)}(\tau; S) \times \mathcal{E}(\sigma; T) \label{eq:bijA2}
\end{align}
\end{thm}

\begin{proof}
We will provide proof for \eqref{eq:bijA1} and remark that the proof of \eqref{eq:bijA2} is nearly identical.

For $\pi \in \mathfrak{S}_n$, we can write the set of all enriched $\pi$-partitions $f: \pi \to S' \times T'$ as follows:
\begin{equation}
\label{eq:epp2}
\begin{aligned}
\mathcal{E}(\pi) = \{\, F =((s_1,t_1),\ldots,(s_n,t_n)) \in (S'\times T')^n & \mid (s_1, t_1) \leq (s_2, t_2) \leq \cdots \leq (s_n, t_n), \\
&  i \notin \Des(\pi) \Rightarrow (s_i, t_i) \leq^{+} (s_{i+1}, t_{i+1})  \\
&  i \in \Des(\pi) \Rightarrow (s_i, t_i) \leq^{-} (s_{i+1}, t_{i+1})  \, \}.
\end{aligned}
\end{equation}
We will now sort the points $F$ into distinct cases. For any $i=1,2,\ldots,n-1$, if $\pi(i) < \pi(i+1)$, then $(s_{i},t_{i}) \leq^{+} (s_{i+1},t_{i+1})$, which falls into one of
two mutually exclusive cases:
\begin{align}
s_i \leq^{+} s_{i+1} & \mbox{ and }  t_{i}\leq^{+} t_{i+1}, \mbox{ or} \label{eqn:e1}\\
s_{i} \leq^{-} s_{i+1} & \mbox{ and }  t_{i} \geq^{-} t_{i+1}. \label{eqn:e2}
\end{align}
If $\pi(i) > \pi(i+1)$, then $(s_{i},t_{i}) \leq^{-} (s_{i+1},t_{i+1})$,
which we split as:
\begin{align}
s_{i} \leq^{+} s_{i+1} & \mbox{ and }  t_{i} \leq^{-} t_{i+1}, \mbox{ or} \label{eqn:e3}\\
s_{i} \leq^{-} s_{i+1} & \mbox{ and }  t_{i} \geq^{+} t_{i+1},\label{eqn:e4}
\end{align}
also mutually exclusive. Define $I_F$ to be the set of all $i$ such that either \eqref{eqn:e2} or \eqref{eqn:e4} holds for $F$. Notice that in both cases, $s_i \leq^{-} s_{i+1}$. Now for any $I \subset [n-1]$, let $A_I$ be the set of all $F$ satisfying $I_F = I$. We have $\mathcal{E}(\pi;S\times T) = \coprod_{I \subset [n-1]} A_I$.

For any particular $I\subset [n-1]$, form the poset $P_{I}$ of the elements
$1,2,\ldots,n$ by $\pi(s) <_{P_{I}} \pi(s+1)$ if $s \notin I$,
$\pi(s) >_{P_{I}} \pi(s+1)$ if $s\in I$. We form a ``zig-zag" poset (see Figure \ref{fig:zigzag}) of $n$ elements labeled consecutively by $\pi(1),
\pi(2),\ldots,\pi(n)$ with downward zigs corresponding to the
elements of $I$.

\begin{figure} [h]
\centering
\includegraphics{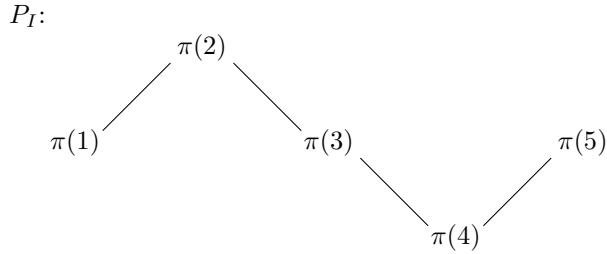}
\caption{The zig-zag poset $P_{I}$ for $I = \{2,3\} \subset
[5]$. \label{fig:zigzag}}
\end{figure}

For any $F$ in $A_I$, let $f: [n] \to T'$ be defined by $f(\pi(i)) = t_i$. It is straightforward to verify that $f$ is an enriched $P_{I}$-partition. Conversely, any enriched $P_{I}$-partition $f$ gives a point $F$ in $A_I$ since by cases \eqref{eqn:e1}--\eqref{eqn:e4} above, if $t_i = f(\pi(i))$, then \[(( s_1,t_1),\ldots,(s_n,t_n)) \in A_I\] if and only if $s_1 \leq \cdots \leq s_n$ and $s_i \leq^{-} s_{i+1}$ for all $i \in I$, $s_i \leq^+ s_{i+1}$ for $i \notin I$. We can therefore turn our attention to enriched $P_{I}$-partitions.

Let $\sigma \in \mathcal{L}(P_{I})$. Recall by Observation \ref{ob1} that $\sigma^{-1}\pi(i) < \sigma^{-1}\pi(i+1)$ if $\pi(i) <_{P_{I}} \pi(i+1)$, i.e., if $i \notin I$. If $\pi(i) >_{P_{I}} \pi(i+1)$ then $\sigma^{-1}\pi(i) > \sigma^{-1}\pi(i+1)$ and $i \in I$. We get that $\Des(\sigma^{-1}\pi) = I$ if and only if $\sigma \in \mathcal{L}(P_I)$. Set $\tau = \sigma^{-1}\pi$. We have \[ \mathcal{E}(\tau; S) = \{ s_{1} \leq \cdots \leq s_{n} \,|\, s_{i} \leq^- s_{i+1} \mbox{ if } i \in \Des(\tau), s_i \leq^+ s_{i+1} \mbox{ otherwise} \}, \] and since $\Des(\tau) = I$, we can write $A_I$ as \[\coprod_{\substack{\sigma\in\mathcal{L}(P_{I}) \\ \sigma\tau = \pi}} \{ F \in (S'\times T')^n \,| \, (s_1, \ldots, s_n) \in \mathcal{E}(\tau;S), (t_{\pi^{-1}\sigma(1)}, \ldots, t_{\pi^{-1}\sigma(n)}) \in \mathcal{E}(\sigma;T) \}. \] Running over all subsets $I \subset [n-1]$, we obtain \[ \mathcal{E}(\pi;S\times T) = \coprod_{\sigma\tau = \pi} \{F \in (S'\times T')^n \,| \, (s_1, \ldots, s_n) \in \mathcal{E}(\tau;S), (t_{\pi^{-1}\sigma(1)}, \ldots, t_{\pi^{-1}\sigma(n)}) \in \mathcal{E}(\sigma;T) \}.\] (Note that $\pi^{-1}\sigma = \tau^{-1}$.) Now we can see the obvious bijection $\coprod_{\sigma\tau=\pi} \mathcal{E}(\tau;S)\times \mathcal{E}(\sigma;T) \to \mathcal{E}(\pi;S\times T)$ given by \[ ( (s_1, \ldots, s_n), (t_1, \ldots, t_n) ) \mapsto ( (s_1, t_{\tau(1)}), \ldots, (s_n, t_{\tau(n)}) ).\]
\end{proof}

\begin{figure} [h]
\centering
\includegraphics[scale=1]{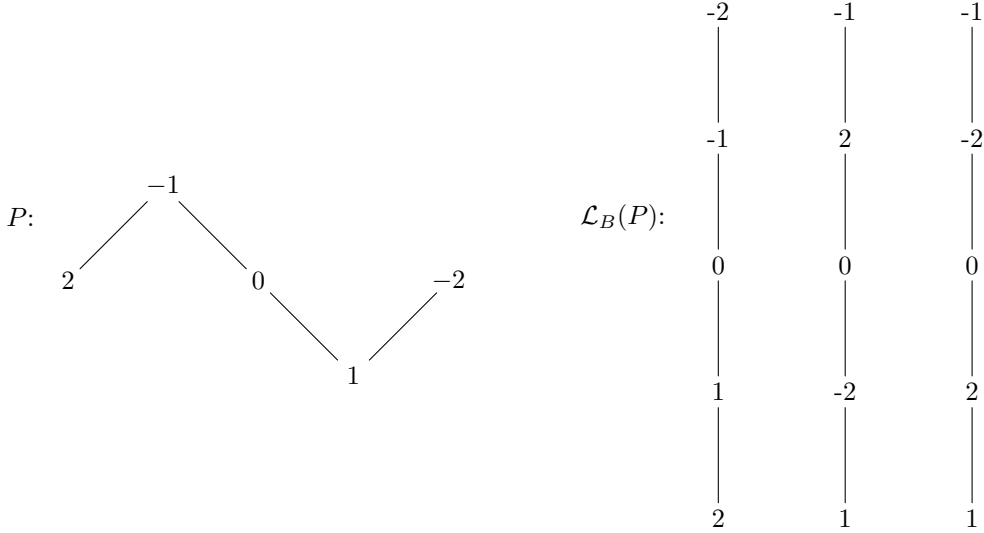}
\caption{A $\mathfrak{B}_{2}$ poset and its linear extensions.\label{fig:typeBlinearextension}}
\end{figure}

Now we present the type B enriched $P$-partitions. When working with signed permutations, we need to change our notion of a poset slightly. See Chow \cite{Chow}; this definition is a simpler version of the notion due to Vic Reiner \cite{Reiner}.
\begin{defn}
A \emph{type B poset}, or \emph{$\mathfrak{B}_{n}$ poset}, is a poset $P$ whose elements
are $0, \pm 1, \pm 2, \ldots, \pm n$ such that if $i <_{P} j$ then
$-j <_{P} -i$.
\end{defn}
Note that if we are given a poset with $n+1$ elements labeled by $0, a_{1},\ldots, a_{n}$ where $a_{i} = i$ or
$-i$, then we can extend it to a $\mathfrak{B}_{n}$ poset of $2n+1$ elements. For example, the $P$ in Figure \ref{fig:typeBlinearextension} could be specified by the relations $0 >_P 1 <_P -2$. In the same way, any signed permutation $\pi \in \mathfrak{B}_{n}$ is a $\mathfrak{B}_{n}$ poset under the total order $\pi(s) <_{\pi} \pi(s+1)$, $0\leq s\leq n-1$. If $P$ is a type B poset, let $\mathcal{L}_B(P)$ denote the set of linear extensions of $P$ that are themselves type B posets. Then $\mathcal{L}_B(P)$ is naturally identified with some set of signed permutations. See Figure \ref{fig:typeBlinearextension}.

We will present some alternate notation for the set $S'$ introduced above. Let $S = \{s_1, s_2, \ldots\}$ be any countable totally ordered set. Then we define the set $S'$ to be the set \[\{ s_1^{-1}, s_1, s_2^{-1}, s_2, \ldots\},\] with total order \[s_1^{-1} < s_1 < s_2^{-1} < s_2 < \cdots \] We introduce this new notation because we want to avoid confusion in defining the set \[\mathbb{Z}' = \{ \ldots, -2, -2^{-1}, -1, -1^{-1}, 0, 1^{-1}, 1, 2^{-1}, 2, \ldots \},\] with the total order \[ \cdots -2 < -2^{-1} < -1 < -1^{-1} < 0 < 1^{-1} < 1 < 2^{-1} < 2 < \cdots\] In general, if we define $\pm S = \{ \ldots, -s_2, -s_1, s_0, s_1, s_2, \ldots\}$, we have the total order on $\pm S'$ given by \[ \cdots -s_{2} < -s_{2}^{-1} < -s_{1} < -s_{1}^{-1} < s_0 < s_{1}^{-1} < s_{1} < s_{2}^{-1} < s_{2} < \cdots\]

For any $s$ in $\pm S'$, let $\varepsilon(s)$ be the exponent on $s$, and let $|s|$ be a map $\pm S' \to S$ that forgets signs and exponents. For example, if $s = -s_{i}^{-1}$, then $\varepsilon(s) = -1 < 0$ and $|s| = s_{i}$, while if $s = s_{i}$, then $\varepsilon(s) = 1 > 0$ and $|s| = s_i$. For $i=0$, we require $\varepsilon(s_0) = 1 > 0$, $|s_0| = s_0$, and $-s_0 = s_0$. Let $s \leq^{+} t$ mean that $s < t$ in $\pm S'$ or $s = t$ and $\varepsilon(s) > 0$. Similarly define $s \leq^{-} t$ to mean that $s < t$ in $\pm S'$ or $s=t$ and $\varepsilon(s) < 0$.

\begin{defn}[Type B enriched $P$-partition]
For any $\mathfrak{B}_n$ poset $P$, an \emph{enriched $P$-partition of type B} is an order-preserving map $f : \pm[n] \to \pm S'$ such that for every $i <_{P} j$ in $P$,
\begin{enumerate}
\item $f(i) \leq^{+} f(j)$ only if $i < j$ in $\mathbb{Z}$,
\item $f(i) \leq^{-} f(j)$ only if $i > j$ in $\mathbb{Z}$,
\item $f(-i) = -f(i)$.
\end{enumerate}
\end{defn}

This definition differs from type A enriched $P$-partitions only in the last condition. It forces $f(0) = s_0$, and if we know where to map $a_1, a_2, \ldots, a_n$, where $a_i = i$ or $-i$, then it tells us where to map everything else. Let $\mathcal{E}_B(P;S)$ denote the set of all type B enriched $P$-partitions $f: P \to \pm S'$. There is a fundamental lemma for type B.

\begin{lem}
We have, \[\mathcal{E}_B(P) = \coprod_{\pi \in \mathcal{L}_B(P)} \mathcal{E}_B(\pi).\]
\end{lem}

We can easily characterize the type B enriched $\pi$-partitions in terms of descent sets, keeping in mind that if we know where to map $i$, then we know where to map $-i$ by the symmetry property: $f(-i) = -f(i)$. For any signed permutation $\pi \in \mathfrak{B}_n$ we have
\begin{equation}\label{eq:eppsetB}
\begin{aligned}
\mathcal{E}_B(\pi) = \{\, f: [n] \to \pm S' & \mid  s_0 \leq f(\pi(1)) \leq f(\pi(2)) \leq \cdots \leq f(\pi(n))\\
 &  \quad i \notin \Des_B(\pi) \Rightarrow f(\pi(i)) \leq^{+} f(\pi(i+1)),  \\
 &  \quad i \in \Des_B(\pi) \Rightarrow f(\pi(i)) \leq^{-} f(\pi(i+1)) \,\}.
\end{aligned}
\end{equation}
Notice that since $\varepsilon(s_0) = 1$, then $s_0 \leq^{-} f(\pi(1))$ is the same as saying $s_0 < f(\pi(1))$, and $s_0 \leq^{+} f(\pi(1))$ is the same as $s_0 \leq f(\pi(1))$. We will show that the set of type B enriched $P$-partitions relates to the set of type B peaks. First, we present the main theorem for type B enriched $P$-partitions. Its proof varies only slightly from that of Theorem \ref{thm:bijA} and is omitted. Let $\mathcal{E}_B(P;S\times T)$ denote the set of all enriched $P$-partitions $f: P \to \pm S' \times \pm T'$ with the up-down order.

\begin{thm}\label{thm:bijB}
We have the following bijection:
\begin{equation}
\mathcal{E}_B(\pi;S\times T) \longleftrightarrow \coprod_{\sigma\tau = \pi} \mathcal{E}_B(\tau; S) \times \mathcal{E}_B(\sigma; T) \label{eq:bijB}
\end{equation}
\end{thm}

\section{Generating functions}\label{sec:qsym}

Recall that a quasisymmetric function is a formal series \[Q(x_1, x_2, \ldots ) \in \mathbb{Z}[[x_1, x_2,\ldots ]] \] of bounded degree such that for any composition $\alpha = (\alpha_1, \alpha_2, \ldots, \alpha_k)$, the coefficient of $x_{1}^{\alpha_1} x_{2}^{\alpha_2} \cdots x_{k}^{\alpha_k}$ is the same as the coefficient of $x_{i_1}^{\alpha_1} x_{i_2}^{\alpha_2} \cdots x_{i_k}^{\alpha_k}$ for all $i_1 < i_2 < \cdots < i_k$. Recall that a composition of $n$, written $\alpha \models n$, is an ordered tuple of positive integers $\alpha = (\alpha_1, \alpha_2, \ldots, \alpha_k)$ such that $|\alpha| = \alpha_1 + \alpha_2 + \cdots + \alpha_k = n$. In this case we say that $\alpha$ has $k$ parts, or $l(\alpha) = k$. We can put a partial order on the set of all compositions of $n$ by refinement. The covering relations are of the form \[ (\alpha_1, \ldots, \alpha_i + \alpha_{i+1}, \ldots, \alpha_k ) \prec (\alpha_1, \ldots, \alpha_i, \alpha_{i+1}, \ldots, \alpha_k).\] Let $\Q_n$ denote the set of all quasisymmetric functions homogeneous of degree $n$. Then $\Q := \bigoplus_{n \geq 0} \Q_n$ denotes the graded ring of all quasisymmetric functions, where $\Q_0 = \mathbb{Z}$.

The most obvious basis for $\Q_n$ is the set of \emph{monomial} quasisymmetric functions, defined for any composition $\alpha = (\alpha_1, \alpha_2, \ldots, \alpha_k) \models n$,
\[ M_{\alpha} := \sum_{i_1 < i_2 < \cdots < i_k} x_{i_1}^{\alpha_1} x_{i_2}^{\alpha_2} \cdots x_{i_k}^{\alpha_k}.\]
There are $2^{n-1}$ compositions of $n$, and hence, the graded component $\Q_n$ has dimension $2^{n-1}$ as a vector space. We can form another natural basis with the \emph{fundamental} quasisymmetric functions, also indexed by compositions,
\[ F_{\alpha} := \sum_{ \alpha \leq \beta } M_{\beta},\] since, by inclusion-exclusion we can express the $M_{\alpha}$ in terms of the $F_{\alpha}$:
\[ M_{\alpha} = \sum_{ \alpha \leq \beta } (-1)^{l(\beta) - l(\alpha)} F_{\beta}.\]
There is a well-known bijection between compositions of $n$ and subsets of $[n-1]$ given by \[ \alpha \mapsto I(\alpha) = \{ \alpha_1, \alpha_1+\alpha_2, \ldots, \alpha_1 + \cdots + \alpha_{k-1} \},\] and so we can also write $M_{\alpha} = M_{I(\alpha)}$ or $F_{\alpha} = F_{I(\alpha)}$ when convenient.

Define the generating function for enriched $P$-partitions $f: P \to \mathbb{P}'$ by \[ \dd(P) = \sum_{f \in \mathcal{E}(P)}\prod_{i=1}^{n} z_{|f(i)|}.\] Then clearly $\dd(P)$ is a quasisymmetric function. By the fundamental Lemma \ref{lem:FLEPP}, we have that \[\dd(P) = \sum_{\pi \in \mathcal{L}(P) } \dd(\pi). \]

For any subset of the integers $I$, define the set $I+1 = \{ i+1 \mid i \in I \}$. From \cite{Stembridge} we see that the generating function for enriched $\pi$-partitions depends only on the peak set of $\pi$.

\begin{thm}[Stembridge \cite{Stembridge}, Proposition 2.2]\label{thm:stem}
For $\pi \in \mathfrak{S}_n$, we have the following equality: \[ \dd(\pi) = \sum_{\substack{ E \subset [n-1]
\\ \Pe(\pi) \subset E \cup (E+1) } } {\kern -10pt} 2^{|E|+1}M_{E}. \]
\end{thm}
For any sets $I$ and $J$, let $I \vartriangle J = (I \cup J)\setminus (I \cap J)$ denote the symmetric difference of sets. The generating functions are also $F$-positive.
\begin{thm}[Stembridge \cite{Stembridge}, Proposition 3.5]\label{thm:peakdesgf}
For $\pi \in \mathfrak{S}_n$, we have the following equality:
\begin{equation}\label{eq:peakdesgf}
 \dd(\pi) = 2^{|\Pe(\pi)|+1} {\kern -10pt} \sum_{\substack{ D \subset [n-1] \\ \Pe(\pi) \subset D \vartriangle (D + 1)}} {\kern -10pt} F_D.
\end{equation}
\end{thm}

For interior peak sets $I$, let $K_I$ be the quasisymmetric function defined by \[ K_{\Pe(\pi)} := \dd(\pi).\] Let $\mathbf{\Pi}_n$ denote the space of quasisymmetric functions spanned by the $K_I$, where $I$ runs over all interior peak sets of $[n-1]$. Stembridge then defines the set of \emph{peak functions} $\mathbf{\Pi} := \bigoplus_{n\geq 0} \mathbf{\Pi}_n$, which is a graded subring of $\Q$. He proved that the functions $K_I$ are linearly independent, and so the rank of $\mathbf{\Pi}_n$ is the the number of distinct interior peak sets, which happens to be the Fibonacci number $f_{n-1}$, defined by $f_0 = f_1 = 1$ and $f_n = f_{n-1} + f_{n-2}$ for $n \geq 2$.

Before discussing generating functions for left enriched $P$-partitions and type B enriched $P$-partitions, we need to introduce Chow's type B quasisymmetric functions \cite{Chow}. Define a \emph{pseudo-composition} of $n$, written $\alpha \Vdash n$, to be an ordered tuple of nonnegative integers $(\alpha_1, \alpha_2, \ldots, \alpha_k)$ whose sum $|\alpha| = \alpha_1 + \cdots + \alpha_k$ is $n$, where $\alpha_1 \geq 0$, $\alpha_i > 0$ for $i > 1$. In other words, given any ordinary composition $\alpha \models n$, we have two corresponding pseudo-compositions: $\alpha$ and $0\alpha = (0, \alpha_1, \ldots, \alpha_k)$. The partial order on the set of all pseudo-compositions of $n$ is again by refinement.

Now we can define a type B quasisymmetric function to be a formal series \[Q(x_0, x_1, x_2, \ldots ) \in \mathbb{Z}[[x_0, x_1, x_2,\ldots ]] \] of bounded degree such that for any pseudo-composition $\alpha = (\alpha_1, \alpha_2, \ldots, \alpha_k)$, the coefficient of $x_{0}^{\alpha_1} x_{1}^{\alpha_2} \cdots x_{k-1}^{\alpha_k}$ is the same as the coefficient of $x_{0}^{\alpha_1} x_{i_2}^{\alpha_2} \cdots x_{i_k}^{\alpha_k}$ for all $0 < i_2 < \cdots < i_k$. Let $\BQ_n$ denote the set of all quasisymmetric functions homogeneous of degree $n$. Then $\BQ := \bigoplus_{n\geq 0} \BQ_n$ is the ring of type B quasisymmetric functions. As before we have a monomial and fundamental basis for $\BQ_n$. For any pseudo-composition $\alpha = (\alpha_1, \alpha_2, \ldots, \alpha_k) \Vdash n$, the monomial functions are
\[
 M_{B,\alpha} := \sum_{ i_2 < \cdots < i_k} x_{0}^{\alpha_1} x_{i_2}^{\alpha_2} \cdots x_{i_k}^{\alpha_k}.\]
There are $2^n$ pseudo-compositions of $n$, so the dimension of $\BQ_n$ is $2^n$. The fundamental basis is
\[ F_{B,\alpha} := \sum_{ \alpha \leq \beta } M_{B,\beta}.\]
There is a bijection between pseudo-compositions of $n$ and subsets of $[0,n-1]$ given by the same map \[ \alpha \mapsto I(\alpha) = \{ \alpha_1, \alpha_1+\alpha_2, \ldots, \alpha_1 + \cdots + \alpha_{k-1} \},\] and so we can also write $M_{B,\alpha} = M_{B,I(\alpha)}$ or $F_{B,\alpha} = F_{B,I(\alpha)}$ when convenient.

Define the generating functions for left enriched $P$-partitions $f: P \to \mathbb{P}^{(\ell)}$, and type B enriched $P$-partitions $f: P \to \mathbb{Z}'$,
\begin{align*}
\dd^{(\ell)}(P) & =  \sum_{f \in \mathcal{E}^{(\ell)}(P)} \prod_{i =1}^{n} z_{|f(i)|},\\
\dd_B(P) & =  \sum_{f \in \mathcal{E}_B(P)} \prod_{i =1}^{n} z_{|f(i)|}.
\end{align*}
The fundamental lemma gives that
\begin{align*}
\dd^{(\ell)}(P) & =  \sum_{\pi \in \mathcal{L}(P)} \dd^{(\ell)}(\pi),\\
\dd_B(P) & =  \sum_{\pi \in \mathcal{L}_B(P)} \dd_B(\pi).
\end{align*}
We can relate $\dd^{(\ell)}(\pi)$ and $\dd_B(\pi)$ to the monomial and fundamental quasisymmetric functions of type B. Notice that for a permutation $\pi \in \mathfrak{S}_n \subset \mathfrak{B}_n$, left peaks coincide with the type B peaks. Therefore we can view left enriched $P$-partitions as a special case of type B enriched $P$-partitions. Furthermore, since Stembridge's enriched $P$-partitions are simply those left enriched $P$-partitions that are nonzero, we have \[ \dd(P)(z_1, z_2, \ldots) = \dd^{(\ell)}(P)(0,z_1,z_2,\ldots ),\] so the results for $\dd(P)$ can be obtained from our results for $\dd^{(\ell)}(P)$ by setting $z_0 = 0$.

\begin{thm}\label{thm:mon}
For $\pi \in \mathfrak{B}_n$, we have the following equations:
\begin{align*}
\dd_B(\pi) & = \sum_{ \substack{ E \subset [0,n-1] \\ \Pe_B(\pi) \subset E \cup (E+1) }} {\kern -10pt} 2^{|E|} M_{B,E}, \\
 & = 2^{|\Pe_B(\pi)|} {\kern -10pt} \sum_{ \substack{ D \subset [0,n-1] \\ \Pe_B(\pi) \subset D \vartriangle (D+1) }} {\kern -10pt} F_{B,D}.
\end{align*}
\end{thm}
We omit the proof of this theorem, but remark that it follows the same lines of reasoning as in Stembridge's proofs of Theorems \ref{thm:stem} and \ref{thm:peakdesgf}.

\begin{cor}
For $\pi \in \mathfrak{S}_n$, we have the following equations:
\begin{align*}
\dd^{(\ell)}(\pi) & = \sum_{ \substack{ E \subset [0,n-1] \\ \lPe(\pi) \subset E \cup (E+1) }} {\kern -10pt} 2^{|E|} M_{B,E},\\
 & = 2^{|\lPe(\pi)|} {\kern -10pt} \sum_{ \substack{ D \subset [0,n-1] \\ \lPe(\pi) \subset D \vartriangle (D+1) }} {\kern -10pt} F_{B,D}.
\end{align*}
\end{cor}

\begin{cor}
The function $\dd_B(\pi)$ depends only on the type B peak set of $\pi \in \mathfrak{B}_n$, the function $\dd^{(\ell)}(\pi)$ depends only on the left peak set of $\pi\in \mathfrak{S}_n$.
\end{cor}

We define the functions $K_{B,I}$ by \[K_{B, \Pe_B(\pi)} := \dd_B(\pi).\]
Note that for a permutation $\pi \in \mathfrak{S}_n$, if $\dd^{(\ell)}(\pi) = K_{B,I}$ then $0 \notin I$.

Let $\mathbf{\Pi}_{B,n}$ denote the span of the $K_{B,I}$, where $I$ ranges over all type B peak sets of $[0,n-1]$. It is not hard to see that the $K_{B,I}$ are linearly independent, and so by counting the number of type B peak sets we see $\mathbf{\Pi}_{B,n}$ has rank $f_{n+1}$. If we define the \emph{type B peak functions}, $\mathbf{\Pi}_B := \bigoplus_{n \geq 0} \mathbf{\Pi}_{B,n}$, then we can see it is a subring of $\BQ$, as an argument identical to that of \cite{Stembridge} Theorem 3.1 shows.

Similarly, let $\mathbf{\Pi}^{(\ell)}_n$ denote the span of all $K_{B,I}$, where $I$ ranges over the left peak sets in $[1,n-1]$. Then $\mathbf{\Pi}^{(\ell)}_n$ has rank $f_n$ and the \emph{left peak functions}, $\mathbf{\Pi}^{(\ell)}:= \bigoplus_{n\geq 0} \mathbf{\Pi}^{(\ell)}_n$, form a subring of $\mathbf{\Pi}_B$.

\section{Duality}\label{sec:coalgebras}

Let $X = \{ x_1, x_2, \ldots \}$ and $Y = \{ y_1, y_2, \ldots \}$ be two sets of commuting indeterminates. Define the set $XY = \{ xy : x \in X, y \in Y\}$. Then we define the bipartite generating function, \[ \dd(P)(XY) = \sum_{ F \in \mathcal{E}(P; \mathbb{P}\times \mathbb{P})} {\kern -10pt} x_{s_1}\cdots x_{s_n} y_{t_1} \cdots y_{t_n}. \] The functions $\dd^{(\ell)}(P)(XY)$ and $\dd_B(P)(XY)$ are defined similarly. Then the following are consequences of Theorem \ref{thm:bijA} and Theorem \ref{thm:bijB}.

\begin{thm}
For any $\pi \in \mathfrak{S}_n$, we have the following equations:
\begin{align}
\dd(\pi)(XY) & = \sum_{\sigma\tau = \pi} \dd(\tau)(X) \dd(\sigma)(Y), \label{eq:internal}\\
\dd^{(\ell)}(\pi)(XY) & = \sum_{\sigma\tau = \pi} \dd^{(\ell)}(\tau)(X) \dd^{(\ell)}(\sigma)(Y). \label{eq:internal2}
\end{align}
\end{thm}

\begin{thm}
For any $\pi \in \mathfrak{B}_n$, we have the following equation:
\begin{equation}\label{eq:internalB}
\dd_B(\pi)(XY) = \sum_{\sigma\tau = \pi} \dd_B(\tau)(X) \dd_B(\sigma)(Y).
\end{equation}
\end{thm}

The formulas above imply duality between $\mathbf{\Pi}_n$ and $\mathfrak{P}_n$, $\mathbf{\Pi}^{(\ell)}_n$ and $\mathfrak{P}^{(\ell)}_n$, and $\mathbf{\Pi}_{B,n}$ and $\mathfrak{P}_{B,n}$. Moreover, they give an explicit combinatorial description for the structure constants of the algebras. We will show how this works for the case of interior peaks. The steps of the construction are the same for the other cases.

First, notice that equation \eqref{eq:internal} implies that \[ K_C (XY) = \sum_{ A,B} c_{A,B}^{C} K_A(X) K_B(Y),\] where the sum is over all pairs of interior peak subsets $A$ and $B$ of $[2,n-1]$, and if $\pi \in \mathfrak{S}_n$ is any permutation with $\Pe(\pi) = C$, then $c_{A,B}^C$ is the number of pairs of permutations $\sigma, \tau$ such that $\Pe(\sigma) = B$, $\Pe(\tau) = A$, and $\sigma\tau = \pi$. We now use this formula to define $\mathbf{\Pi}_n$ as a coalgebra with coproduct $\Delta: \mathbf{\Pi}_n \to \mathbf{\Pi}_n \otimes \mathbf{\Pi}_n$ defined as \[ \Delta( K_C ) = \sum_{ A, B} c_{A,B}^C K_A \otimes K_B.\] We can define a coalgebra $\mathbb{Z}[\mathfrak{S}_n]^*$ dual to the group algebra with coproduct defined as \[ \Delta(\pi) = \sum_{\sigma\tau = \pi} \tau \otimes \sigma.\] Define the map $\varphi^* : \mathbb{Z}[\mathfrak{S}_n]^* \to \mathbf{\Pi}_n$ by $ \varphi^*(\pi) = K_{\Pe(\pi)}$, which, by \eqref{eq:internal}, is a surjective homomorphism of coalgebras. Now we dualize.

Let $\mathbf{\Pi}_n^*$ be the algebra dual to $\mathbf{\Pi}_n$, with basis elements $K_I^*$. By definition, multiplication in this basis is \[ K_A^* K_B^* = \sum_{ C} c_{A,B}^C K_C^*,\] where the sum is over all interior peak subsets $C$. The dual of $\varphi^*$ is now an injective homomorphism of algebras, $\varphi: \mathbf{\Pi}_n^* \to \mathbb{Z}[\mathfrak{S}_n]$ defined by \[ \varphi(K_I^*) = \sum_{ \Pe(\pi) = I} \pi = v_I.\] Thus the interior peak algebra can be defined as the image of $\varphi$, and the structure constants carry through: \[ v_A v_B = \sum_{C} c_{A,B}^C v_C.\]

We describe the structure constants for the left and type B peak algebras:
\begin{itemize}
\item Let $d_{A,B}^C$, over triples of left peak sets $A,B,C$, be the structure constants for $\mathfrak{P}_n^{(\ell)}$. Then for any $\pi \in \mathfrak{S}_n$ such that $\lPe(\pi) = C$, $d_{A,B}^C$ is the number of pairs of permutations $\sigma, \tau$ such that $\lPe(\sigma) = B$, $\lPe(\tau) = A$, and $\sigma\tau = \pi$.

\item Let $e_{A,B}^C$, over triples of type B peak sets $A,B,C$, be the structure constants for $\mathfrak{P}_{B,n}$. Then for any $\pi \in \mathfrak{B}_n$ such that $\Pe_B(\pi) = C$, $e_{A,B}^C$ is the number of pairs of permutations $\sigma, \tau$ such that $\Pe_B(\sigma) = B$, $\Pe_B(\tau) = A$, and $\sigma\tau = \pi$.
\end{itemize}

The type B peak algebra is something new, defined as the linear span of sums of signed permutations in $\mathfrak{B}_n$ with common type B peak set.
\begin{thm}
The space $\mathfrak{P}_{B,n}$ is a subalgebra of $\mathbb{Z}[\mathfrak{B}_n]$ of dimension $f_{n+1}$. (In fact it is a subalgebra of Solomon's descent algebra for type $B_n$.)
\end{thm}

\begin{remark}
While $\mathfrak{P}^{(\ell)}_n$ was introduced in \cite{AguiarBergeronNyman}, the authors had no combinatorial description for its structure constants (they were working over a different basis), and neither were the structure constants for $\mathfrak{P}_n$ known. Independently, Nantel Bergeron and Christophe Hohlweg \cite{BergeronHohlweg} recently found the same description we give here.
\end{remark}

\begin{remark} Theorem \ref{thm:bijA} can be modified to combine left enriched $P$-partitions and interior enriched $P$-partitions. When translated to generating functions, it implies that $\mathbf{\Pi}_n$ is a two-sided ideal in $\mathbf{\Pi}_n^{(\ell)}$, and hence $\mathfrak{P}_n$ is an ideal in $\mathfrak{P}^{(\ell)}_n$.
\end{remark}

\section{Specializations}\label{sec:spec}

Define the polynomial $\Omega(P;x)$, called the \emph{enriched order polynomial}, by \[\Omega(P; k) := \dd(P)(\underbrace{1,1,\ldots,1}_{k}, 0,0,\ldots),\] meaning we set $z_i = 1$ for $i = 1,\ldots, k$, and $z_i = 0$ for $i > k$. It turns out that for $\pi \in \mathfrak{S}_n$, $\Omega(\pi; x)$ is a polynomial of degree $n$ that only depends on the \emph{number} of interior peaks of $\pi$. We can use order polynomials to study commutative peak algebras, spanned by sums of permutations with the same \emph{number} of peaks. We sketch the idea for the interior peaks case.

Let $E_i$ be the sum of all permutations with $i$ interior peaks, and let $\Omega(i;x)$ denote the order polynomial for any such permutation.  Now we define:
\begin{align*}
 \rho(x) & :=  \sum_{\pi \in \mathfrak{S}_n} \Omega(\pi;x/2) \pi =  \sum_{i=1}^{\lfloor \frac{n+1}{2} \rfloor} \Omega(i;x/2)E_{i} =
\begin{cases}
\displaystyle \sum_{i=1}^{n/2} e_i x^{2i} & \mbox{ if $n$ is even, }\\
\displaystyle \sum_{i=1}^{(n+1)/2} e_i x^{2i-1} & \mbox{ if $n$ is odd. }
\end{cases}
\end{align*}
The function $\rho(x)$ is a polynomial in $x$ with coefficients in the group algebra $\mathbb{Q}[\mathfrak{S}_n]$ (we now need to work over the rational numbers). From Theorem \ref{thm:bijA} we can obtain the following.

\begin{thm}
As polynomials in $x$ and $y$ with coefficients in the group algebra $\mathbb{Q}[\mathfrak{S}_n]$, we have:
\[\rho(x)\rho(y)  = \rho(xy).\]
\end{thm}

What the theorem tells us is that the coefficients $e_i$ are mutually orthogonal idempotents. With a little more work, we see that the span of the $E_i$ is the same as the span of the $e_i$, so that the sums of permutations with common peak numbers span a commutative $\lfloor \frac{n+1}{2} \rfloor$-dimensional subalgebra of the group algebra. This algebra and its left peak variant were introduced in \cite{AguiarBergeronNyman}.

\begin{table}[h]
\begin{tabular}{c| c c c c}
($\delta_{ij}e^{*}_i e^{*}_j$) & $e_j$ & $\overline{e}_j$ & $e^{(\ell)}_j$ & $e^{(r)}_j$ \\
\hline

$e_i$ &  $e_i$ & $e_i$ & $e_i$ & $e_i$ \\

$\overline{e}_i$  & $\overline{e}_i$ & $\overline{e}_i$ & $\overline{e}_i$ & $\overline{e}_i$ \\

$e^{(\ell)}_i$ & $e_i$ & $\overline{e}_i$ & $e^{(\ell)}_i$ & $e^{(r)}_i$ \\

$e^{(r)}_i$ & $\overline{e}_i$ & $e_i$ & $e^{(r)}_i$ & $e^{(\ell)}_i$

\end{tabular}
\vspace{.5cm}
\caption{Multiplication table for type A coefficients. \label{table:2}}
\end{table}

We can use the same approach to get similar results for other commutative peak algebras, given by the span of sums of permutations with the same number of right peaks, exterior peaks, and type B peaks. Table \ref{table:2} summarizes how the different type A peak idempotents interact (though $e_i^{(r)}$ is not technically idempotent). We remark that while the number of right peaks does not give a basis on its own, its multiplicative closure is still a proper subalgebra.

\section{Negative results}\label{sec:neg}

Before anything was proved, the type B peak algebra was found experimentally, and along the road to its discovery there were several dead-end definitions. To save others the trouble of these detours, we finish with a list of some subalgebras that \emph{do not} exist in general. For the symmetric group, the sums of permutations with the same set of right peaks do not form an algebra, nor do the sums of permutations with the same set of exterior peaks. For the hyperoctahedral group, we do not get a proper subalgebra by taking the sums of permutations with the same: interior peak set (i.e., ignoring peaks at 0), number of interior peaks, exterior peak set, or exterior peak number.

\end{document}